\pdfoutput=1
\documentclass[a4paper, reqno, twocolumn]{amsart}

\usepackage{graphicx}
\usepackage{amssymb}
\usepackage{amsmath}
\usepackage{hyperref}
\long\def\symbolfootnote[#1]#2{\begingroup%
\def\thefootnote{\fnsymbol{footnote}}\footnote[#1]{#2}\endgroup}
\begin{document}

%\footnotetext[1]{Research supported by the National Science Foundation through the Research Experiences for Undergraduates Program at Cornell.}
%\footnotetext[2]{Research supported in part by the National Science Foundation, grant DMS-0652440.}

\title{Localized Eigenfunctions: here you see them, there you don't}
\author{Steven M. Heilman}

\author{Robert S. Strichartz}

\thanks{\textit{Steven M. Heilman is a graduate student of mathematics at the Courant Institute of Mathematical Sciences, New York University.  He was supported by the National Science Foundation through the Research Experiences for Undergraduates Program at Cornell.  His email address is} \texttt{heilman@cims.nyu.edu}\\\\
\textit{Robert S. Strichartz is a professor of mathematics at Cornell University.  He was supported in part by the National Science Foundation, grant DMS-0652440.  His email address is} \texttt{str@math.cornell.edu}}

\maketitle

%%%%%%%%%%%%%%%%%%%%%%%%%%%%%
% Address and Email Section %
%%%%%%%%%%%%%%%%%%%%%%%%%%%%%

%-Describe Eigenfunction Localization
%-Show first localized eigenfunction
%    (Look further up in spectrum?)
%    (Could be example of Non-QUE domain, but Hassel's results may not be directly adaptable...)
%
%-Show power laws for parameter variation
%
%-Show other symmetric domain, display same results

The Laplacian $\Delta$ ($\frac{\partial^{2}}{\partial x^{2}}+\frac{\partial^{2}}{\partial y^{2}}$ in the plane) is one of the most basic operators in all of mathematical analysis.  It can be used to construct the important spacetime equations of mathematical physics, such as the heat equation, the wave equation and the Schr\"{o}dinger equation of Quantum Mechanics.  It has been studied from many points of view, and in many different contexts (Riemannian manifolds, graphs, fractals, etc).  The eigenvalue equation\\ $-\Delta u=\lambda u$ in a domain $\Omega$ with suitable boundary conditions (Dirichlet conditions $u|_{\partial\Omega}=0$ and Neumann conditions $\frac{\partial u}{\partial n}|_{\partial\Omega}=0$ are the most common) defines both the eigenvalue $\lambda$ (in many cases a nonnegative real number) and the eigenfunction $u$.  The set of all eigenvalues (usually a discrete set) is called the \textit{spectrum} of the Laplacian on $\Omega$, and there is a vast literature on the relationship between the spectrum and the geometry of $\Omega$.  See for example \cite{kac71}, \cite{gordon92} and \cite{zelditch09}. The spectrum contains a lot of information about the domain, and teasing out this information, including relationships between quantum and classical mechanics on the domain, has been a fascinating, ongoing and highly nontrivial enterprise.  For example, the famous Weyl law gives the asymptotic size of the $k$th eigenvalue as $k\to\infty$.  There are also interesting and intricate estimates on the size of the first and second nontrivial eigenvalues.  Needless to say, there are no connections between estimates of small and large eigenvalues.

Eigenfunctions are also fascinating objects.\\  Clearly these are more complicated objects, so it can be expected that it is more difficult to say things about them.  Many of the classical special functions and orthogonal polynomials can be related to eigenfunctions of the Laplacian for specific geometries.  There have been many recent studies concerning eigenfunctions associated to large eigenvalues.  In this note we invite you to look at some startling pictures of some specific ``localized'' eigenfunctions associated to small eigenvalues.

\section*{Examples of Localization}
How should we define a ``localized'' eigenfunction?  We would be tempted to say it is an eigenfunction with support $\Omega_{1}$ that is considerably smaller than all of $\Omega$.  But it is well-known that eigenfunctions are real analytic functions, hence cannot vanish on any open set.  So we must be content with saying that the function is ``very small'' on the complement of $\Omega_{1}$.  This is of course not a mathematical definition, although it might be acceptable to a physicist, or a Justice of the Supreme Court.  One could make it into a precise definition with a parameter $\epsilon$ to quantify the statement ``very small'', but this just begs the question: how small does $\epsilon$ have to be to make the statement interesting?  In this note we will show you some pictures to try to convince you that there are surprising examples where $\epsilon$ is smaller than you might expect.

Localized eigenfunctions have been observed before.  As usual, physicists know more than mathematicians in the subject, but with less certainty [this is a kind of uncertainty principle].  Regardless, the cross-pollination in this subject between these two groups over the past century merits recognition and esteem. For high frequency eigenfunctions, relations between eigenfunction localization and billiard dynamics have been studied.  A nonexhaustive list includes treatments of (non)localization on: closed stable geodesics \cite{babich67} and closed unstable geodesics \cite{heller84,verdiere94,burq05,hassell09}.  Other results address dichotomies \cite{berry77}, numerical aspects \cite{barnett07,backer03}, rarity \cite{shnirelman74,verdiere85,zelditch87} and near-nonexistence \cite{lindenstrauss06} of such (phase space) localized eigenfunctions, as $\lambda\to\infty$.  In the low frequency realm, no deep explanation for eigenfunction localization seems to exist.  Low frequency, or ``ground state'' eigenfunctions have been widely studied (see for example \cite{payne73,banuelos00}.  However, the authors can only find scattered examples of low frequency localization, such as: near a fractal boundary \cite{sapoval97}, in narrow channels between domains \cite{courant34} and in square pairs with irrational ratios of frequency oscillations (Example $3$ in \cite{jones08}).  As is well known, an eigenfunction is itself an eigenfunction in each of its nodal domains (with appropriate boundary conditions).  Therefore, results for low frequency eigenfunctions can inform higher frequency studies.

    \begin{figure}[htbp!]
    \begin{center}
    \includegraphics[scale=.4]{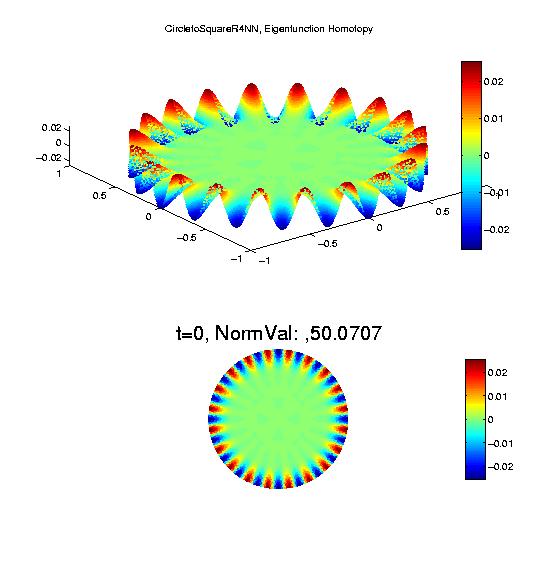}
    \end{center}
    \caption{Localized Circle Eigenfunction}
    \label{figzero}
    \end{figure}

A simple example of an eigenfunction on the disc localized to a neighborhood of the boundary circle is shown in Figure \ref{figzero}.  These examples tend to involve eigenfunctions with eigenvalues rather high up in the spectrum and domains with special types of billiard flows.  In contrast, our examples occur low down in the spectrum and are consequences of symmetry considerations.  We work with Neumann boundary conditions because they are natural (the weak formulation of the eigenvalue equation is $-\int_{\Omega}(\nabla u\cdot\nabla v)\,dx=\lambda\int_{\Omega}uv\,dx$ for all reasonable test functions $v$, without imposing any boundary conditions) and they were essential in our work on approximating fractal Laplacians with ordinary planar Laplacians \cite{berry08}, \cite{heilman09}.  Similar examples with Dirichlet boundary conditions also exist.  It was our coauthor Tyrus Berry who first serendipitously discovered examples of localized eigenfunctions on sawtooth shaped domains as reported in \cite{berry08}, but these examples did not play any role in the theory developed there.  By coincidence, many localized eigenfunctions of fractal Laplacians have been known since the work of Fukushima and Shima \cite{fukushima92}, and these can also be explained by symmetry considerations \cite{barlow97}.  See \cite{strichartz99,strichartz06} for expository accounts of this phenomenon.

    \begin{figure}[htbp!]
    \begin{center}
    \includegraphics[scale=.6]{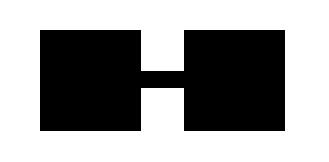}
    \end{center}
    \caption{Two rooms with a passage}
    \label{fig.5}
    \end{figure}

    \begin{figure}[htbp!]
    \begin{center}
    \includegraphics[scale=.35]{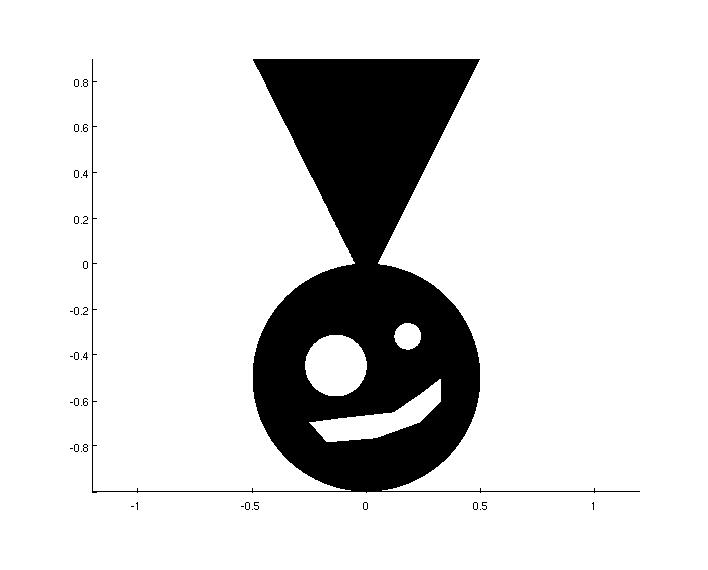}
    \end{center}
    \caption{Smiley Domain, height $h=0.1$}
    \label{figone}
    %Eigenvalue ?
    \end{figure}

    \begin{figure}[htbp!]
    \begin{center}
    \includegraphics[scale=.36]{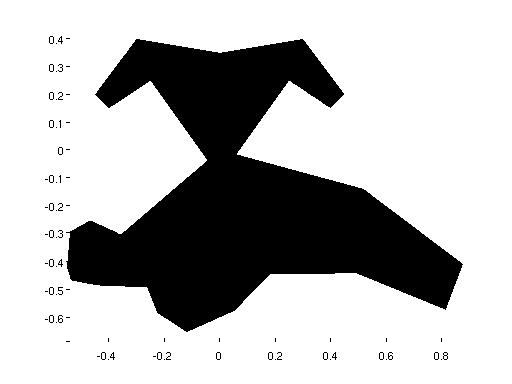}
    \end{center}
    \caption{Cow Domain, height $h=0.1$}
    \label{figseven}
    %Eigenvalue ?
    \end{figure}

Our examples can be though of as modified versions of ``rooms and passages'' domains \cite{courant34}.  If $\Omega=\Omega_{1}\cup\Omega_{2}\cup\Omega_{3}$ consists of two rooms, $\Omega_{1},\Omega_{2}$, and a very short narrow connection passage $\Omega_{3}$ (see Figure \ref{fig.5}) then it is not surprising that there are Dirichlet eigenfunctions on $\Omega$ that are very close to Dirichlet eigenfunctions on $\Omega_{1}$ extended to be zero on $\Omega_{2}\cup\Omega_{3}$.  In our examples we will join two rooms $\Omega=\Omega_{1}\cup\Omega_{2}$ that intersect in a small, but not very small piece.  Figure \ref{figone} shows one example, the ``smiley face'', and Figure \ref{figseven} shows another, the ``cow''.  The main idea is that $\Omega_{1}$ must possess an axis of symmetry, and the rooms join together at one end of this axis.  In other words, there is a line $L$ such that the reflection $R$ in $L$ preserves $\Omega_{1}$, $R\Omega_{1}=\Omega_{1}$, and $L$ passes through $\Omega_{1}\cap\Omega_{2}$.  Suppose $u_{1}$ is a Neumann eigenfunction of the Laplacian on $\Omega_{1}$ that is skew-symmetric with respect to $R$, $u_{1}(Rx)=-u_{1}(x)$.  Such eigenfunctions occur throughout the spectrum, since indeed every eigenspace splits into functions that have symmetry and skew-symmetry with respect to $R$ (most eigenspaces are $1$-dimensional and are of one or the other symmetry type).  Then $u_{1}$ vanishes along $L$.  If the point $p$ where $L$ intersects $\partial\Omega_{1}$ is a corner point, then $u_{1}$ and $\nabla u_{1}$ vanish at $p$, so $u_{1}$ is relatively small near $p$.  (For example, $u_{1}(x,y)=\cos\pi x-\cos\pi y$ at the origin if $\Omega_{1}$ is the unit square.)  So it is not surprising that there is an eigenfunction $u$ on $\Omega$ that is close to $u_{1}$ on $\Omega_{1}$ and close to zero off $\Omega_{1}$.

\begin{figure}[htbp!]
    \begin{center}
    \includegraphics[scale=.5]{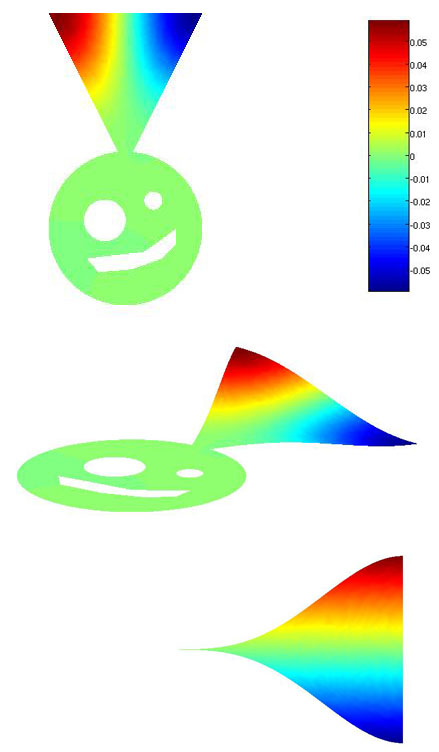}
    \end{center}
    \caption{Smiley, Fifth Eigenfunction}
    \label{figtwo}
    %Eigenvalue 5
    \end{figure}

\begin{figure}[htbp!]
    \begin{center}
    \includegraphics[scale=.25]{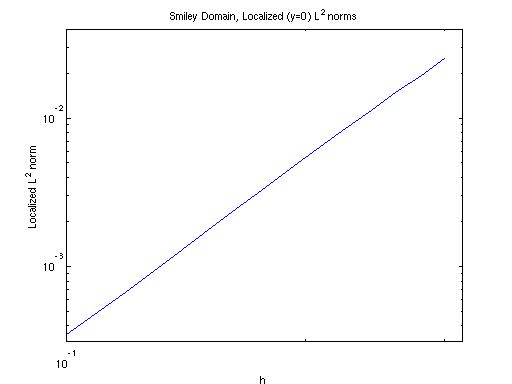}
    \end{center}
    \caption{Localization in $L^{2}$ norm, for Figure \ref{figtwo} example}
    \label{figfour}
    \end{figure}

\begin{figure}[htbp!]
    \begin{center}
    \includegraphics[scale=.25]{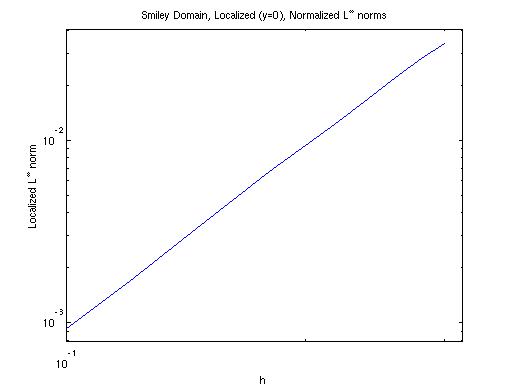}
    \end{center}
    \caption{Localization in $L^{\infty}$ norm, for Figure \ref{figtwo} example}
    \label{figfive}
    \end{figure}

\begin{figure}[htbp!]
    \begin{center}
    \includegraphics[scale=.5]{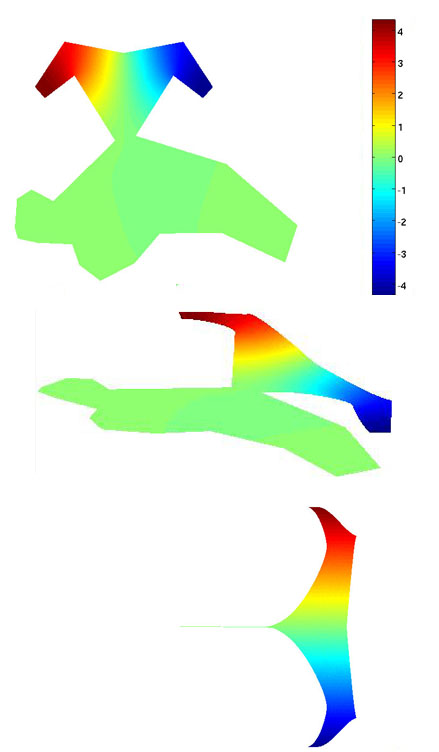}
    \end{center}
    \caption{Cow, Fourth Eigenfunction}
    \label{figeight}
    %Eigenvalue ?
    \end{figure}

\begin{figure}[htbp!]
    \begin{center}
    \includegraphics[scale=.25]{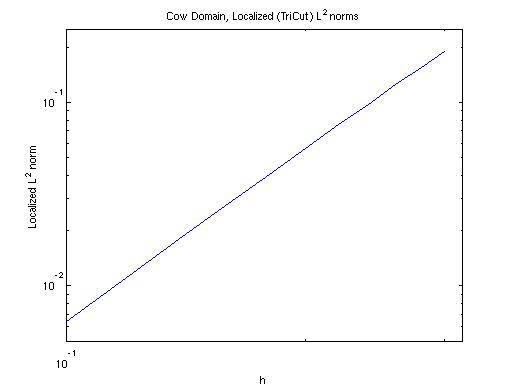}
    \end{center}
    \caption{Localization in $L^{2}$ norm, for Figure \ref{fignine} example}
    \label{fignine}
    \end{figure}

\begin{figure}[htbp!]
    \begin{center}
    \includegraphics[scale=.25]{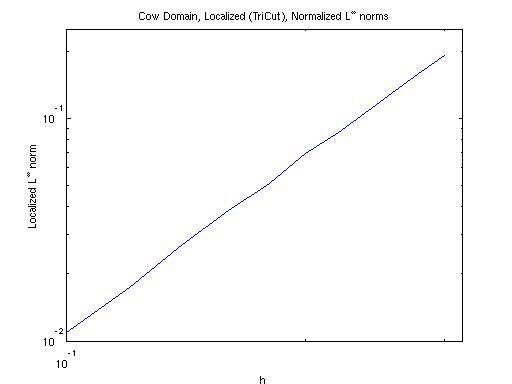}
    \end{center}
    \caption{Localization in $L^{\infty}$ norm, for Figure \ref{fignine} example}
    \label{figten}
    \end{figure}

\begin{figure}[htbp!]
    \begin{center}
    \includegraphics[scale=.45]{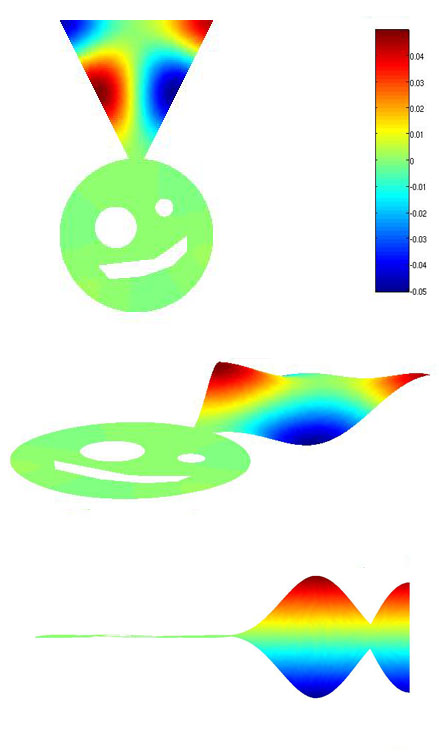}
    \end{center}
    \caption{Smiley, Twelfth Eigenfunction}
    \label{figeleven}
    %Eigenvalue ?
    \end{figure}

\begin{figure}[htbp!]
    \begin{center}
    \includegraphics[scale=.45]{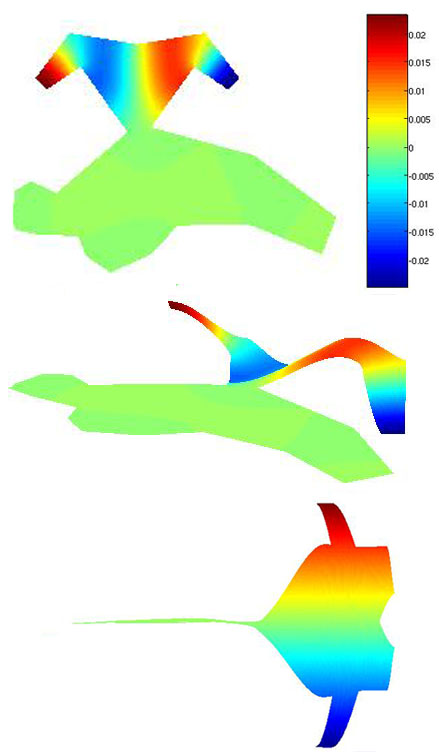}
    \end{center}
    \caption{Cow, Eleventh Eigenfunction}
    \label{figtwelve}
    %Eigenvalue ?
    \end{figure}

Such reasoning does not yield a sharp estimate for how localized $u$ is, so we look at some experimental evidence.  We use Matlab to numerically approximate some eigenfunctions on our domains using the finite element method.  Figures \ref{figtwo}-\ref{figtwelve} show the results.  Note that we normalize the eigenfunctions to have $L^{2}$ norm on $\Omega$ equal to $1$, and we can measure the localization either by the $L^{2}$ norm on $\Omega\setminus\Omega_{1}$, a kind of average localization, or the $L^{\infty}$ norm on $\Omega\setminus\Omega_{1}$, a uniform localization.  In Figures \ref{figfour},\ref{figfive},\ref{fignine}, and \ref{figten} we show both of these as a function of the aperture size [height.. will edit later] for the connection for the each domain.  These are log-log plots, suggesting a power law relationship over the given range of $h$ values.  In the tables below we give the best fit power law for two eigenfunctions in each domain. Note that the powers vary considerably in these four examples.
%We even find a power law relationship for the rotational angle of the triangle.
\begin{table}[htbp!]
\resizebox{7cm}{!}{
\begin{tabular}{c|cc}
 \multicolumn{3}{p{6cm}}{Summary Table: Smiley Domain}\\
\hline
 & Eigfcn $5$ & Eigfcn $12$\\
\hline
$L^{2}$ Localization      & $y=11.254x^{3.9087}$ & $y=249.06x^{2.4636}$\\
$L^{\infty}$ Localization & $y=4.1735x^{3.2959}$ & $y=90.552x^{2.3553}$\\
\hline
\end{tabular}
}
%\caption{
%Summary Table
%}
\label{tableone}
\end{table}

\begin{table}[htbp!]
\resizebox{7cm}{!}{
\begin{tabular}{c|cc}
 \multicolumn{3}{p{5cm}}{Summary Table: Cow Domain}\\
\hline
 & Eigfcn $4$ & Eigfcn $11$\\
\hline
$L^{2}$ Localization      & $y=119.65x^{3.0889}$ & $y=2096.5x^{2.8028}$\\
$L^{\infty}$ Localization & $y=31.615x^{2.6075}$ & $y=676.08x^{2.5700}$\\
\hline
\end{tabular}
}
%\caption{
%Note: Eigfcn $11$ plots not so linear
%}
\label{tabletwo}
\end{table}

\section*{Conclusion}
Our examples show how easy it is to find surprisingly localized eigenfunctions.  The domains do not have to have any special properties beyond the symmetry of the $\Omega_{1}$ piece (breaking the symmetry even slightly makes the examples disappear).  We do not have to go very high up in the spectrum.  As mathematicians it is natural for us to want a theorem that explains the examples, or at least a conjectured theorem.  Perhaps there is such a theorem, and a perceptive reader might be able to find one, but at present we don't see any.  Or, we might suggest that there is more to mathematics than just theorems.  This might sound like a radical suggestion, or perhaps it is just common sense.

%%%%%%%%%%%%%%%%
% Bibliography %
%%%%%%%%%%%%%%%%

%We would also like to keep the figures and tables above the references so I am inserting the following command
\clearpage
%\bibpunct{[}{]}{,}{b}{,}{,}

%%%%%%%%%%%%%%%%%%
% Figure Section %
%%%%%%%%%%%%%%%%%%

\end{document}